\documentclass[12pt]{article}
\usepackage{amsmath,amsthm,amsfonts,amssymb}
\date{}

\newcommand{\ind}{\mathrm{ind}}
\newcommand{\IR}{\mathbb R}

\newtheorem{problem}{Problem}

\begin{document}

\title{Some Open Problems in Topological Algebra}

\author{Taras Banakh, Mitrofan Choban, Igor Guran, Igor Protasov}

\maketitle

This is the list of open problems in topological algebra posed on
the conference dedicated to the 20th anniversary of the Chair of
Algebra and Topology of Lviv National University, that was held on
28 September 2001.
\bigskip

\begin{problem}[Choban] Is every topological group a quotient group of a
zero-dimensional topological group of the same weight? When an
almost metrizable topological group  is  a quotient of a
zero-dimensional group of the same weight?
\end{problem}

A space $X$ is zero-dimensional if $\dim X=0$. A topological group
is {\em almost metrizable} if contains a compact subset of
countable character. Let us make some comments to this problem. If
a topological group $G$ is a quotient of an almost metrizable
group $H$ with $\ind H=0$, then $G$ contains a zero-dimensional
compact subgroup of countable character in $G$.

If $H$ is a zero-dimensional subgroup of countable character in a
group $G$ and $xh=hx$ for every $x\in G$ and $h\in H$, then $G$ is
a quotient group of some almost metrizable zero-dimensional group
of the same weight. In particular the answer to the first part of
Problem 1 is positive for metrizable groups (see \cite{Ch2},
\cite{Be}) and for almost metrizable abelian groups.
A.V.Arhangel'skii proved that every topological group is a
quotient of a zero-dimensional $\sigma$-discrete group. For
universal algebras this fact was proved in \cite{Ch1}. \smallskip

\begin{problem}[Choban] Under which conditions is the free
universal algebra of an uncountable signature over a metrizable
space $X$ paracompact? In particular, is the free topological
linear space $L(X)$ of $X$ over the discrete field of real numbers
paracompact?
\end{problem}
\bigskip

A semigroup $S$ with the identity $e$ endowed with a topology is
called a {\em left} (resp. {\em right}) {\em bounded} if for every
neighborhood $U$ of $e$ there is a finite subset $F$ of $S$ such
that $S=FU$ (resp. $S=UF)$; $S$ is called {\em bounded} if $S$ is
both left and right bounded.

For a topological space $X$ let $S(X)\subset X^X$ be the semigroup
of all continuous selfmappings of $X$ endowed with the topology of
pointwise convergence (i.e., the topology inherited from the
Tychonov product $X^X$). A topological space $X$ is called {\em
homogeneous} if for any points $x_1,x_2\in X$ there is a
homemorphism $h$ of $X$ with $h(x_1)=x_2$.

\begin{problem}[Protasov] Is the semigroup $S(X)$ left bounded for every
zero-dimensional compact homogeneous space?
\end{problem}

The answer is positive provided $X$ has a base of the topology,
consisting of pairwise homeomorphic clopen subsets, see
\cite{Pr1}.

\begin{problem}[Protasov] Is the semigroup $S(X)$ right bounded
for any zero-dimensional homogeneous space $X$?
\end{problem}

The answer is positive provided $X$ has a base of the topology,
consisting of clopen subsets homeomorphic to $X$, see \cite{Pr1}.

Let $X$ be a topological space. A subgroup $H$ of $S(X)$ is called
{\em distal} if for any distinct points $x_1,x_2\in X$ and any
point $x\in X$ there is a neighborhood $U$ of $x$ such that
$\{h(x_1),h(x_2)\}\not\subset U$ for all $h\in H$. It is proven in
\cite[Theorem 4]{Pr1} that a left bounded subgroup $H\subset S(X)$
is distal provided $H$ acts transitively on $X$.

\begin{problem}[Protasov] Let $X$ be a compact space and $H$ be a distal
subgroup of $S(X)$, acting transitively on $X$. Is $H$ left
bounded?
\end{problem}

Under a {\em left-topological group} we understand a pair
$(G,\tau)$ consisting of a group $G$ and a topology $\tau$
invariant with respect to the left shifts $l_g:x\mapsto gx$. If,
in addition $\tau$ is invariant with respect to the right shifts,
then $(G,\tau)$ is called a {\em semitopological group}. A
semitopological group $G$ with continuous inverse mapping
$x\mapsto x^{-1}$ is called a {\em quasitopological group}. If the
group operation of $G$ is continuous with respect to the topology
$\tau$, then $(G,\tau)$ is a {\em paratopological group}. If,
additionally, the operation of taking the inverse is continuous,
then $(G,\tau)$ is a topological group.

It is well known that $\sigma$-compact topological groups have
countable cellularity \cite{Tk} while compact topological groups
support a strictly positive probability measure (i.e., a Borel
probability measure $\mu$ such that $\mu(U)>0$ for any nonempty
open subset $U$ of the group). Recently T.~Banakh and O.~Ravsky
\cite{BR} (see also \cite{BR2}) proved that any bounded
paratopological group $G$ has countable cellularity (moreover,
each cardinal of uncountable cofinality is a precaliber of $G$).
On the other hand, according to \cite{Pr2} for every infinite
cardinal $\tau$ there is a left bounded left topological group of
cellularity $\tau$.

\begin{problem}[Protasov]\label{pr5} Let $G$ be a bounded semitopological
(quasitopological) group. Has $G$ countable cellularity?
\end{problem}

To answer Problem~\ref{pr5} in negative it suffices to give a
positive answer to

\begin{problem}[Protasov] Does every zero-dimensional compact
homogeneous space admit the structure of a left topological group?\footnote{This problem was answered in negative by Banakh in [T.Banakh, {\em A homogeneous first-countable zero-dimensional compactum failing to be a left-topological group}, Mat. Stud.29 (2008) 215--217].}
In particular, does the countable power $K^\omega$ of a
first-countable zero-dimensional compact space $K$ admit the
structure of a left topological group?
\end{problem}

Let us remark that according to Motorov Theorem the countable
power $K^\omega$ of each first-countable zero-dimensional compact
space $K$ is homogeneous, see \cite{Ar2}. 

On the other hand, a positive answer to Problem~\ref{pr5} would
follow from the positive answer of

\begin{problem}[Protasov] Does every compact left topological
group $G$ support a strictly positive probability Borel measure?
\end{problem}

There exists a compact left-topological group admitting no
invariant probability Borel measure, see \cite{MP}.

\begin{problem}[Protasov] Can every non-discrete topological group $G$ be algebraically
generated by a nowhere dense subset?\footnote{The answer to Problem 6 is ``No'': under CH the topological group $\IR^{\omega_1}$ contains a dense non-separable subgroup $G$ such that each meager subset of $G$ generates a metrizable separable subgroup of $G$; see Proposition 13.8 in [I.Protasov, T.Banakh. Ball stuctures and colorings of graphs and groups // (Mat. Stud. Monograph Series. 11), VNTL Publ. 2003, 148p.]}
\end{problem}

Let us mention that each countable topological group is
algebraically generated by some closed discrete subset (see
reference in \cite{Pr3}) while every left topological group is
algebraically generated by some subset with empty interior
\cite{Pr3}.
\bigskip

A paratopological group $G$ is {\em left $\omega$-bounded} (resp.
{\em right $\omega$-bounded}) if for any neighborhood $U\subset G$
of the origin there is a countable subset $F\subset G$ with $G=FU$
(resp. $G=UF$); $G$ is {\em $\omega$-bounded} if it is both left
and right $\omega$-bounded.

\begin{problem}[Guran] Is a paratopological group $G$ left
$\omega$-bounded if it is right $\omega$-bounded?
\end{problem}

For topological groups the answers to the last two questions are
in positive, see \cite{Gu}.
\bigskip

A subset $A$ of a topological group $G$ is called {\em
$o$-bounded} if for any sequence $(U_n)_{n\in\omega}$ of
neighborhoods of the origin of $G$ there is a sequence
$(F_n)_{n\in\omega}$ of finite subsets of $G$ such that
$A\subset\bigcup_{n\in\omega}F_nU_n$. It is clear that each
$\sigma$-compact topological group is $o$-bounded while each
$o$-bounded group is $\omega$-bounded. Next, given a subset $A$ of
a topological group $G$, consider the following game OF(A)
(abbreviated from Open-Finite). Two players, I and II, choose at
every step $k\in\omega$ a neighborhood $U_n\subset G$ of the
origin, and a finite subset $F_n$ of $G$, respectively. At the end
of the game, the player II is declared the winner if
$A\subset\bigcup_{n\in\omega}F_nU_n$. It is easy to see that for a
$\sigma$-compact group $G$ the player II has a winning strategy in
the game OF($G$). On the other hand, if a topological group $G$ is
not $o$-bounded, then the player I has a winning strategy in
OF($G$).

\begin{problem}[Banakh] Is there a (metrizable) $o$-bounded
topological group $G$ such that the player I has a winning
strategy in the game OF($G$)?
\end{problem}

Such a group, if exists, cannot be analytic and abelian (more
generally, a SIN-group). We remind that a topological space $X$ is
analytic if it is a metrizable continuous image of a separable
complete metric space. On the other hand, the group ${\mathbb
Z}^\omega$ contains a dense $o$-bounded $G_\delta$-subset $A$ such
that the first player has a winning strategy in the game OF($A$),
see \cite{BZ}.

\begin{problem}[Banakh] Let $A$ be a compact subset of the real
line $\IR$ such that $A-A=\{a-a':a,a'\in A\}$ is a neighborhood of
zero in $\IR$. Has the sum $A+A+A+A$  non-empty
interior in $\IR$?
\end{problem}

\end{document}